\documentclass[
reprint,
superscriptaddress,
amsmath,amssymb,
aps,
prl,
floatfix,
 lengthcheck
]{revtex4-1}
\usepackage[utf8]{inputenc}
\usepackage[T1]{fontenc}
\usepackage[version=4]{mhchem}
\usepackage{graphicx}
\usepackage{placeins}
\usepackage{listings}
\usepackage[usenames, dvipsnames]{xcolor}
\usepackage{url}
\usepackage{hyperref}
\usepackage{ulem}
\usepackage{comment}
\date{\today}
\usepackage{algorithm} 
\usepackage{algpseudocode} 
\definecolor{codegreen}{rgb}{0,0.6,0}
\definecolor{codegray}{rgb}{0.5,0.5,0.5}
\definecolor{codepurple}{rgb}{0.58,0,0.82}
\definecolor{backcolour}{rgb}{0.95,0.95,0.92}

\begin{document}
\author{Anton Vladyka}
\email{anton.vladyka@utu.fi}
\affiliation{University of Turku, Department of Physics and Astronomy, FI-20014 Turun yliopisto, Finland}
\author{Eemeli A. Eronen}
\affiliation{University of Turku, Department of Physics and Astronomy, FI-20014 Turun yliopisto, Finland}
\author{Johannes Niskanen}
\email{johannes.niskanen@utu.fi}
\affiliation{University of Turku, Department of Physics and Astronomy, FI-20014 Turun yliopisto, Finland}
\title{Implementation of the Emulator-based Component Analysis}
\begin{abstract}
We present a PyTorch-powered implementation of the emulator-based component analysis used for ill-posed numerical non-linear inverse problems, where an approximate emulator for the forward problem is known. This emulator may be a numerical model, an interpolating function, or a fitting function such as a neural network. With the help of the emulator and a data set, the method seeks dimensionality reduction by projection in the variable space so that maximal variance of the target (response) values of the data is covered. The obtained basis set for projection in the variable space defines a subspace of the greatest response for the outcome of the forward problem. The method allows for the reconstruction of the coordinates in this subspace for an approximate solution to the inverse problem. We present an example of using the code provided as a Python class.
\end{abstract}
\maketitle

\section{Introduction}
Inverse problems consist of reconstructing input data $\mathbf{X}$ or parameters of a model from a later observation of a dependent response $\mathbf{Y}$. In addition to their rich heritage from pure mathematics \cite{Kaipio2005,Calvetti2018, Uhlmann2014}, many applied inverse problems originate from geophysics \cite{Tarantola1982}. Other modern examples of the applied problems include tomography \cite{Arridge2009} (as a part of the more general image/signal processing problem) and spectroscopy. Inverse problems can be identified in a wealth of research contexts, and their approximate solutions are traditionally searched for by regularization or by Bayesian procedures \cite{Kaipio2005,Mohammad-Djafari2022}.
\par
For example, given an atomistic structure with a charge state, the resulting electronic spectrum (ultraviolet/visible, X-ray, etc.) is defined by the quantum-mechanical electronic Hamiltonian of the system. The formed map is certainly nonlinear in terms of structure, and not necessarily an injective function. Moreover, some structural degrees of freedom may account for huge spectral responses while some may result little or virtually none. As often the case with inverse problems, full reconstruction of structural parameters is an ill-posed problem, and already identification of the relevant structural degrees of freedom may be a worthwhile scientific result. It may likewise be important to filter out structural degrees of freedom without spectral response to achieve a more well-posed inverse problem of reconstruction of structural information from the respective spectra \cite{Vladyka2023}.
\par
Machine learning (ML) opens new possibilities for solving inverse problems. With the backpropagation algorithm \cite{Rumelhart1986} and graphics processing units (GPU), training of extensive neural networks (NN) has become feasible. Because the NN is a universal approximator \cite{Hornik1989,Cybenko1989,Leshno1993}, mimicking many forward problems is possible given that enough data is available for the task. The complete solution of a potentially ill-defined inverse problem may still be difficult, but an approximative one may be obtained by first identifying the reconstructable information. To do this, the fast evaluation of the forward problem by NN may be utilized in iterative procedures. We propose the emulator-based component analysis (ECA) \cite{Niskanen2022} as a tool in finding a fast approximate solution for nonlinear inverse problems. The method has been applied to spectroscopic problems from gaseous \cite{Niskanen2022}, to condensed amorphous \cite{Vladyka2023} and to liquid \cite{Eronen2024} phases using the original implementation, and to a molecular liquid using the notably faster presented implementation \cite{Eronen2024b}.
\par
\begin{figure}[t]
\includegraphics[width=\linewidth]{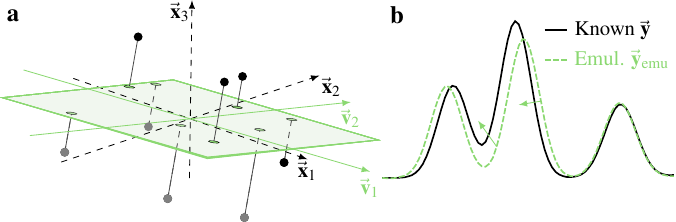}
\caption{The principle of ECA illustrated for 2 ECA components $\{\mathbf{v}_1, \mathbf{v}_2\}$ in a 3-dimensional $\mathbf{X}$ space $\mathrm{span} \{\mathbf{x}_1, \mathbf{x}_2,\mathbf{x}_3\}$ (a). The ECA components (basis vectors) are selected by iterative optimization so that  emulated $\mathbf{y}_\mathrm{emu}$ for the projected data points match closest with the known $\mathbf{y}$ of the original data (b). \label{fig:eca} }
\end{figure}
\par
The ECA algorithm aims for dimensionality reduction by projection on basis vectors in the input space $\mathbf{X}$ while maximizing the variance of a dependent target variable of the forward problem in space $\mathbf{Y}$ of vectors $\mathbf{y}=\mathbf{y}(\mathbf{x})$. Thus ECA belongs to the family of `projection pursuit' algorithms \cite{Kruskal1969,Friedman1974,Huber1985} and is closely related to projection pursuit regression (PPR) \cite{Friedman1981}. However, instead of using an unknown smoothing function (ridge function), ECA utilizes a more general and pre-selected function: a machine learning model ({\it i.e.} a NN) acting as an emulator for the forward problem $ \mathbf{y}_\mathrm{emu}(\mathbf{x})\approx \mathbf{y}(\mathbf{x})$. For example, in spectroscopy studies, such emulator would predict a spectrum for the corresponding atomistic structure. Approximate solution $\mathbf{x}'$ of the inverse problem $\mathbf{y}\rightarrow \mathbf{x}$ may then be found for the obtained dimensionally-reduced space. In this contribution we present an efficient numerical implementation of the ECA algorithm and its application for the inverse problem together with an example.
\par
 
\section{Definition}
The result of ECA is an orthogonal set of vectors $\mathbf{V}=\{\mathbf{v}_j\}_{j=1}^k$ in X space. These vectors are optimized using data $\{(\mathbf{x}_i,\mathbf{y}_i)\}$ separate to those used for emulator training, so that most variance of vectors $\mathbf{y}_i$ is covered by emulator prediction $\mathbf{y}_\mathrm{emu}(\tilde{\mathbf{x}}_i^{(k)})$ subsequent to projection of vectors $\mathbf{x}_i$ onto $\mathbf{V}$
\begin{equation}
\label{projection}
\tilde{\mathbf{x}}_i^{(k)}=\sum_{j=1}^{k}\underbrace{(\mathbf{v}_j\cdot\tilde{\mathbf{x}}_i)}_{=:t_{i,j}}\mathbf{v}_j.
\end{equation}
The algorithm is designed to work with z-standardized $\mathbf{X}$ data points $\tilde{\mathbf{x}}_i$ (corresponding to $\mathbf{x}_i$) for which Equation (\ref{projection}) represents a reduction to $k$ dimensions. The coordinates of the data point $\tilde{\mathbf{x}}^{(k)}_i$ reduced dimension ($\mathbf{t}$ scores) read $\mathbf{t}_i=(t_{i,1},...,t_{i,k})$. The idea is depicted in Figure \ref{fig:eca}.
\par
As a metric for the optimization of $\mathbf{V}$, we use a generalized covered variance (R$^2$-score). Organizing z-standardized Y data points $\tilde{\mathbf{y}}_i$ as row vectors of matrix $\mathbf{Y}$, and the corresponding predicted data points $\mathbf{y}_\mathrm{emu}(\tilde{\mathbf{x}}_i^{(k)})$ as row vectors of matrix $\mathbf{Y}^\mathrm{(pred)}$, this is defined as:
\begin{equation}
\rho = 1 - \frac{\mathrm{tr}\left(\bar{\mathbf{Y}}^T\bar{\mathbf{Y}}\right)}{\mathrm{tr}\left(\mathbf{Y}^T\mathbf{Y}\right)}
\label{R2}
\end{equation}
where $\bar{\mathbf{Y}} = \mathbf{Y} - \mathbf{Y}^\mathrm{(pred)}$. 
\par
The procedure utilizes iterative optimization of the basis vectors so that $\rho$ is maximized for the used data. The method proceeds one vector $\mathbf{v}_j$ at a time where orthonormality is enforced after optimization of each vector. We note that the replacing any of the ECA basis vector $\mathbf{v}_i$ with its opposite $-\mathbf{v}_i$ does not affect the covered variance, but only changes the sign of the corresponding $\mathbf{t}_i$ score. 

\par
The expansion $\tilde{\mathbf{x}}_i^{(k)}$ gives an approximation of point $\tilde{\mathbf{x}}_i$ (and $\mathbf{x}_i$ after inverse z-standardization), in the subspace $\mathrm{span}(\mathbf{V})$, which $\mathbf{Y}$ by design shows greatest response to. We approach the inverse problem by deducing the respective coordinates $\mathbf{t}'$ for a given $\mathbf{y}$, and $\tilde{\mathbf{y}}$, after which the approximate solution is
\begin{equation}
\tilde{\mathbf{x}}'=\sum_{j=1}^k t_j'\mathbf{v}_j.
\end{equation}
This reconstruction utilizes numerical optimization for $\mathbf{t}'$ for given $\mathbf{y}$ and $\mathbf{V}$. The approximate solution $\mathbf{x}'$ in $\mathbf{X}$ space is obtained via z-score inverse transformation.

\section{Implementation}
The module is written in Python utilizing PyTorch \cite{pytorch} (>2.0). The use of PyTorch offers some benefits:
\begin{enumerate}
    \item Runs very efficiently owing to automatic differentiation approach implemented in this package.
    \item Can easily be run in parallel over multiple CPU cores or GPU.
    \item It provides more control during the fitting procedure than some other optimization toolkits, (\textit{e.g.}, scipy.optimize \cite{Scipy}).
\end{enumerate}

We arranged the search for ECA vectors as a mini-batch optimization using the Adam optimizer \cite{Adam}. The algorithm utilizes previously found component vectors, and orthogonality of the search is forced by allowing optimization steps only in their orthogonal complement. The workflow of the procedure is presented in Algorithm \ref{pseudocode} and the optimization options together with their default values are given in Table~\ref{tab:options}. The full reference manual and a sample code for the implementation are available in the Appendix.
\begin{algorithm}[H]
	\begin{algorithmic}
        \Require Trained emulator $\mathbf{y}_\mathrm{emu}: \mathbf{x} \to \mathbf{y} $
        \Require Data $(\mathbf{X, Y})$
        \State Initialize first vector: $\mathbf{v}_1 \gets \mathrm{rnd()}$
        \While {$\mathbf{v}_1$ not optimized}
            \State Calculate projection: $\mathbf{X}_\mathrm{proj} \gets (\mathbf{X}\cdot \mathbf{v}_1)\cdot \mathbf{v}_1$
            \State Optimize $\mathbf{v}_1$ : $\rho(\mathbf{Y}, \mathbf{y}_\mathrm{emu}(\mathbf{X}_\mathrm{proj})) \to max$
            \State Normalize $\mathbf{v}_1 \gets \mathbf{v}_1/|\mathbf{v}_1|$
        \EndWhile
        \State Initialize the basis set: $\mathbf{V}_1 \gets \mathbf{v}_1$
        \For {rank $i \gets 2,\ldots, k$}
            \State Initialize $i$-th component: $\mathbf{v}_i \gets \mathrm{rnd()}$
            \While {$\mathbf{v}_i$ not optimized}
                \State Calculate projection: $\mathbf{X}_\mathrm{proj} \gets (\mathbf{X}\cdot \mathbf{V})\cdot \mathbf{V} + (\mathbf{X}\cdot \mathbf{v}_i)\cdot \mathbf{v}_i$
                \State Optimize $\mathbf{v}_i$ : $\left(\rho(\mathbf{Y}, \mathbf{y}_\mathrm{emu}(\mathbf{X}_\mathrm{proj})) \to max \right)$ $\cap \left(\mathbf{V}\cdot\mathbf{v}_i = 0\right)$  \Comment{$\mathbf{v}_i$ is orthogonal to vector space $\mathbf{V}$}
                \State Normalize $\mathbf{v}_i \gets \mathbf{v}_i/|\mathbf{v}_i|$
            \EndWhile
            \State Append optimized component to the basis: $\mathbf{V}_i \gets \mathbf{v}_i$
        \EndFor
    \end{algorithmic}
\caption{ECA implementation \label{pseudocode}}
\end{algorithm}

\FloatBarrier
\par

\par
Reproducibility of the results is satisfied by setting the seed for the random numbers generator via \texttt{seed} parameter of the options or by calling \texttt{set\_seed()} method. The \texttt{fit()} method also allows the manual step by step evaluation  and re-evaluation of the components by specifying \texttt{keep=n} argument: a positive value is interpreted as `start after $n$-th component', and for negative values the last $|n|$ calculated components are removed, evaluation starts after the remaining ones.

\par
\begin{figure*}
\includegraphics{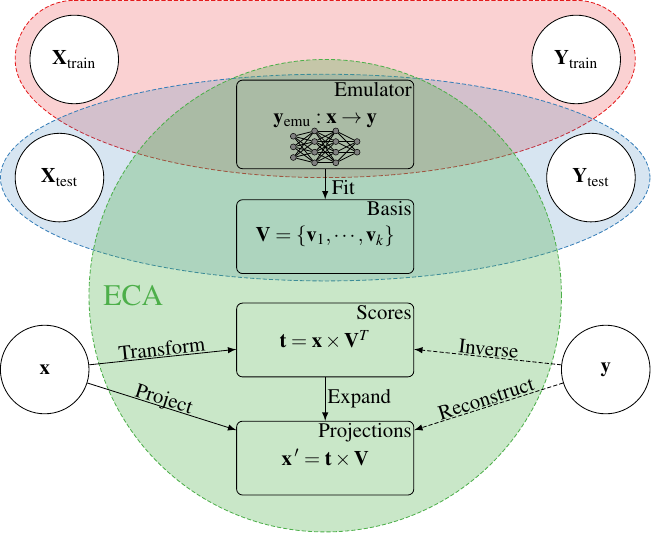}
\caption{Schematics of ECA functionality. ECA uses an emulator, trained on $(\mathbf{X}, \mathbf{Y})_\mathrm{train}$, validated on $(\mathbf{X}, \mathbf{Y})_\mathrm{test}$. Fitting of the ECA object on $(\mathbf{X}, \mathbf{Y})_\mathrm{test}$ (or $(\mathbf{X}, \mathbf{Y})_\mathrm{train}$) results in a basis set $\mathbf{V}$. `Transform' and `project' operations are linear transformations of the input vector $\mathbf{x}$. On the contrary, `inverse' and `reconstruct' are implemented using optimization algorithms. These are \textit{not} linear operations and \textit{may} produce ambiguous results for the complex (highly multidimensional) systems. \label{fig:schematic}}
\end{figure*}

In addition to providing the ECA basis vectors for the data--emulator pair, the presented implementation carries out typically used transformations. These functionalities are presented in Figure~\ref{fig:schematic}. 

\begin{table*}[]
    \caption{The list of options for both forward and inverse algorithms. \label{tab:options}}
    \centering
    \begin{tabular}{l l l l}
    Parameter\hspace{1cm}             & Explanation\hspace{5cm}  & Default\hspace{1cm} & (inv)\\
    \hline
    \hline
        lr            & Learning rate for the Adam optimizer& $10^{-3}$ & $0.05$\\
        betas         & Beta parameters for the Adam optimizer (see \cite{Adam}) & (0.9, 0.999) & (0.9, 0.999) \\
        tol           & Stopping condition & $10^{-4}$ & $10^{-4}$\\
        epochs        & Max. number of epochs & $10^4$ & $10^3$\\
        batch\_size   & Mini-batch size (forward only) & 200 & ---\\
        seed          & Seed for random number generator & None & None\\
    \hline
\end{tabular}
\end{table*}
In this work, we propose to utilize dimensionality reduction achieved via projection of the $\mathbf{X}$ data onto ECA space for approximation of the inverse problem solution. The \texttt{inverse()} method is  searching for the points in ECA space ($\mathbf{t}$ scores) for which the emulator provides the best match with given $\mathbf{y}$, while the \texttt{reconstruct()} method is performing inverse transformation followed by its expansion to $\mathbf{X}$ space. In other words, it returns an $\mathbf{x}'$, which predicts to the best match with given $\mathbf{y}$. The obtained $\mathbf{t}'$ scores can be used for studying the differences in the data, \textit{e.g.}, via possible grouping, cluster analysis etc, and interpretability of $\mathbf{x}'$ vectors allows to find the relation between $\mathbf{x}$ and $\mathbf{y}$ in a significantly more lightweight way than an exact solution of the inverse problem would.

The PyTorch implementation, while robust and efficient, requires an emulator implemented in PyTorch. Numerous machine learning libraries exist and provide a similar functionality. A simple feed-forward network can be transferred between the implementations by specifying the matrices, bias vectors and activation functions. We provide a script which allows to convert a trained sklearn-based multilayer perceptron (\texttt{sklearn.MLPRegressor} \cite{scikit-learn}) into an equivalent PyTorch-based neural network for use with the PyTorch-based ECA.

\section*{Rudimentary Example}
ECA is intended to find a few degrees of freedom that potentially cover most of target variance of the forward problem. As the first example for this we study the problem for the forward function
\begin{equation}
y = y(\mathbf{x}) = (\mathbf{v}\cdot\mathbf{x})^3\label{eq:rudimentary1}
\end{equation}
tailored to represent the problem well. Function $y$ of Equation (\ref{eq:rudimentary1}) is dependent on a single degree of freedom in the input space, and in this example we use
\begin{equation}
\mathbf{v} =\frac{1}{\sqrt{d}}\sum_{i=1}^d \mathbf{\hat{e}}_i ,\label{eq:rudimentary2}
\end{equation}
where $d$ is the dimensionality of the input vectors for the forward problem and $\mathbf{\hat{e}}_i$ the $i$th basis vector of this space. This function depends only on vector $\mathbf{v}$ which expected as the result for first ECA vector $\mathbf{v}_1$. Equally obviously, at most the scores $t_1=\mathbf{v}_1\cdot\mathbf{x}$ are reconstructable quantities for the inverse problem for this non-bijective map, which manifests partial loss of information about $\mathbf{x}$ in the target variable $y$. We do all procedures following the case of an intended study.
\par
We use 20\,000 data points $(\mathbf{x}_i, y_i)$ and vary the dimension $d$ of the input vector space. For simplicity we sample the $\mathbf{x}$ data from a normal distribution, allowing us to skip the element-wise z-score standardization of the $\mathbf{x}_i$ data, otherwise expected at this point. To bring the target values to the scale of typical application, we z-score standardized the $y$ data. We split the data to 80\% for training the emulator (16\,000 points) and 20\% (4\,000 points) for testing and ECA fitting. For the emulator required by ECA, we chose an NN with four hidden layers of 16 neurons and ReLU as the activation function. We implemented the network in PyTorch and trained it using Adam as the optimizer with weight decay of 0.001 and initial learning rate of 0.001. The use of this fixed architecture limits the changes of evaluation time to the ECA procedure and to the inevitable differences in the input-layer operations. For each input-dimensionality $d$ we trained the model with 80\% of the training data and determined the condition for patient early stopping with the rest. In the end R$^2$ score for the test set was always at least 0.94. 
\par
To study the dependency of the ECA fit time on input dimensionality $d$, we ran the fit algorithm 25 times for each $d$ with varying initial guesses. We repeated this procedure for both the presented and for the original implementation using identical emulators. The original implementation uses MLPRegressor class of the sklearn \cite{scikit-learn} package, to which we converted the trained PyTorch emulators. We ran all computations using a single core on a cluster computer (Puhti, CSC, Kajaani, Finland) to minimize hardware-induced variation in the run time. The limited number of trials (25) is due to the slow evaluation time of the original implementation at large $d$. 
\par
\begin{figure}
    \centering
    \includegraphics[width=\columnwidth]{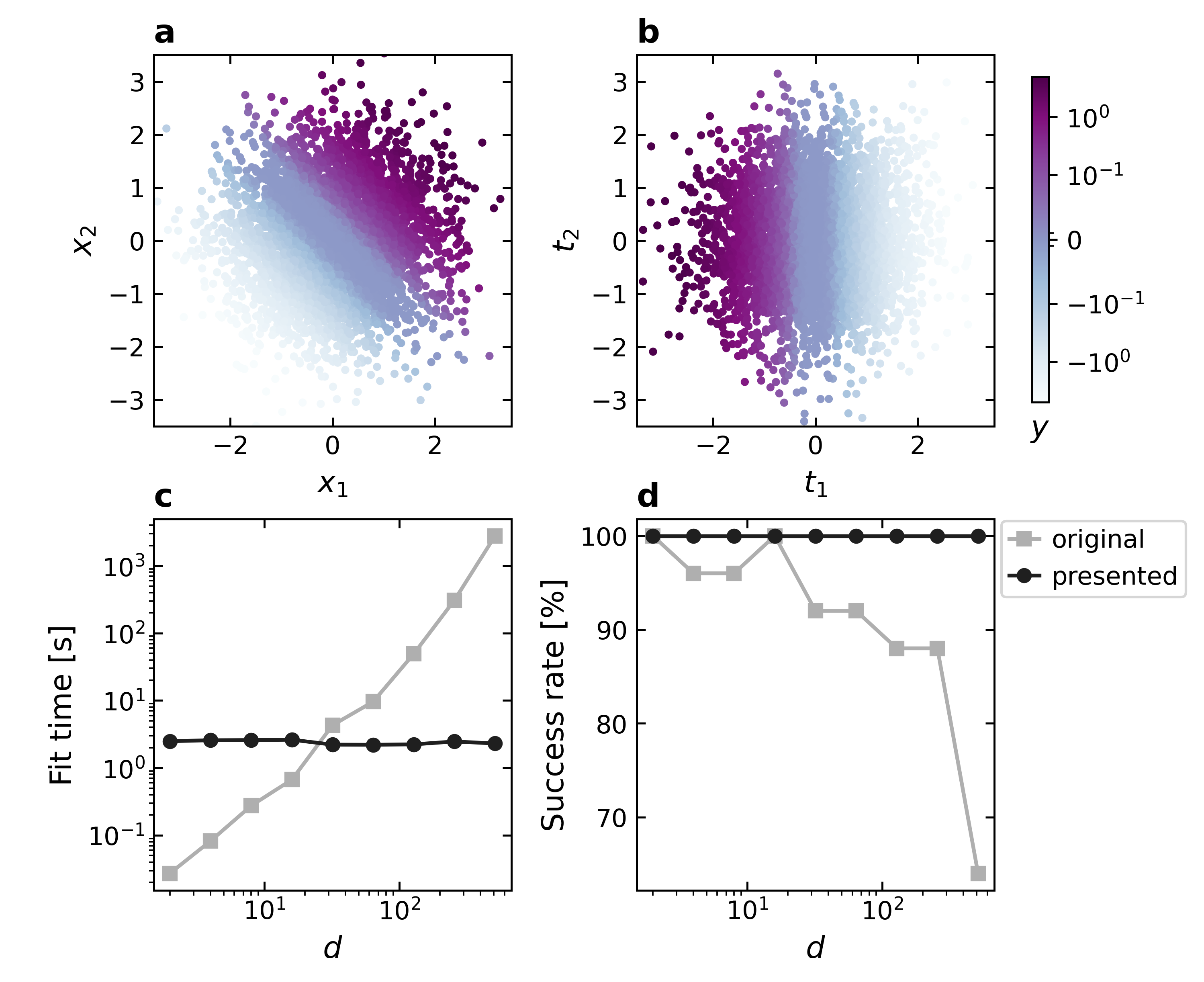}
    \caption{\label{fig:forward} Performance of the new implementation on the rudimentary example data. (a)~Test set $\mathbf{x}=(x_1,x_2)$ for two-dimensional $(d=2)$ case. Color of each point depicts the corresponding value $y(\mathbf{x})$ (Eq. (\ref{eq:rudimentary1})) (b)~Projection scores $\mathbf{t}=(t_1,t_2)$ of the data on the two ECA vectors using the same color scheme. (c)~Evaluation times for ECA fit for $\mathbf{v}_1$ as a function of dimensionality $d$. We observed a linear dependence of the ECA fitting time on the number of data points for both implementations. (d)~Fit success rate measured as fraction of fits with $|\mathbf{v}\cdot\mathbf{v}_1|>0.95$ as a function of $d$ from 25 independent trials.}
\end{figure}

Figure \ref{fig:forward}a represents the test set data points in the case $d=2$, while figure \ref{fig:forward}b depicts this data projected on the two ECA vectors. In this case the vector $\mathbf{v}_2$ is given by the orthonormality condition apart from its sign, but in a higher-dimensional case the choice would not be unique due to equal irrelevance of all other degrees of freedom of the input. Figure \ref{fig:forward}c depicts the evaluation time for the ECA fit with the two implementations as a function of the input dimension $d$. The original implementation based on optimization by SciPy \cite{Scipy} optimization routines are faster in fewer dimensions but the mini-batch optimization routine of the presented implementation shows no significant dependence on $d$, surpassing the former at around $d=30$. We observed a linear dependence of the ECA fitting time on the number of fit data points for both implementations. Figure \ref{fig:forward}d depicts the rate of success of the procedure as a function of $d$, measured as fraction of fits with $|\mathbf{v}\cdot\mathbf{v}_1|>0.95$. The presented implementation always fits $\mathbf{v}_1$ successfully while the original implementation is less stable failing with a third of the initial guesses at $d = 512$. We note that when dimensionality $d$ approaches the number of used data points (in this case 4\,000 points), the fit became increasingly unstable for the presented implementation as well.
\par

\begin{figure}
\includegraphics[width=\columnwidth]{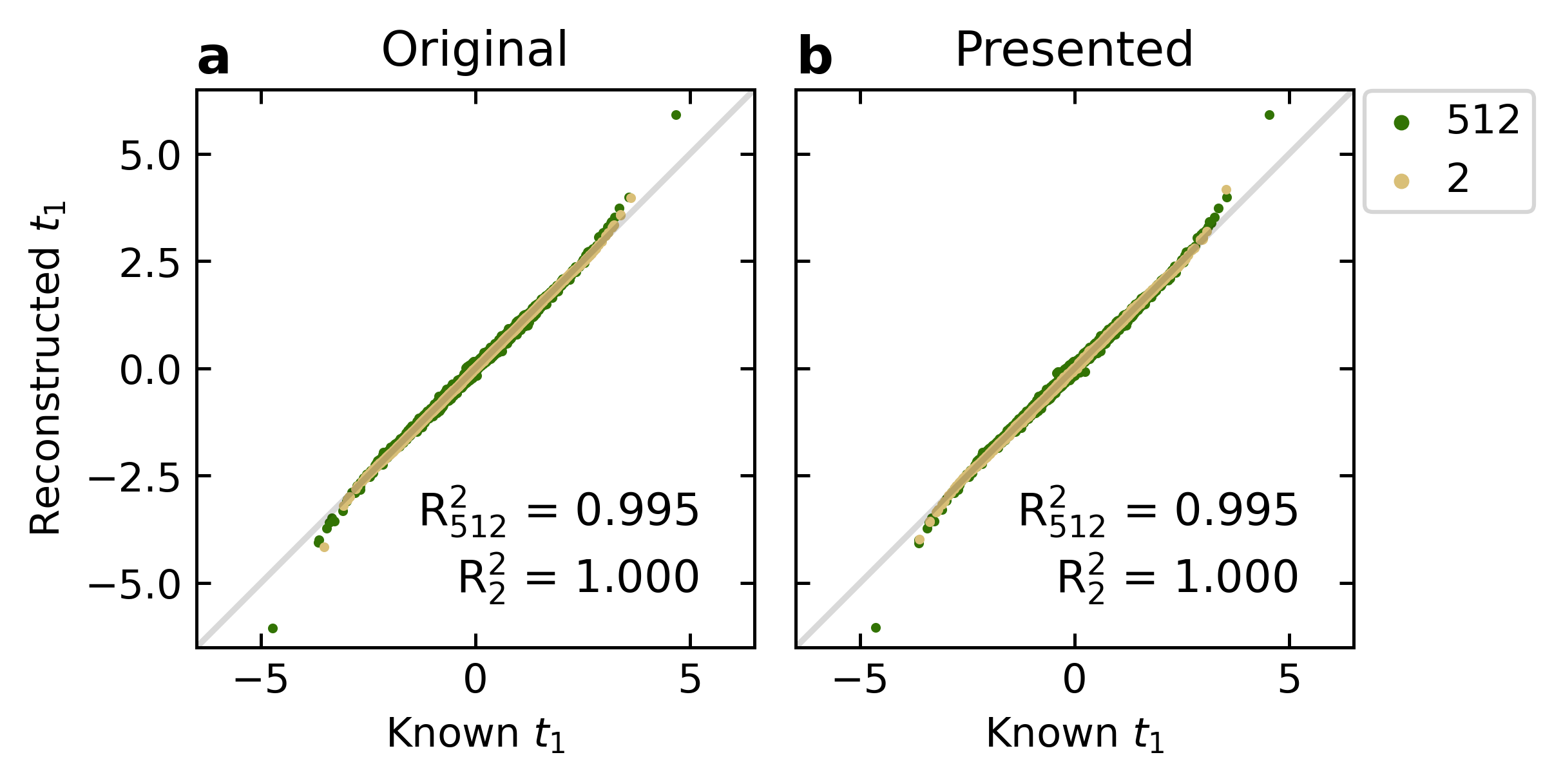}
\caption{\label{fig:inverse}Reconstructed $t_1$ scores versus the known (projected) ones for input dimensionalities $d=2$ and $d=512$ obtained using (a) the original and (b) the presented implementation. The diagonal grey line indicates a perfect match. The scatter plots show that both the implementations work well in solving the inverse problem for $t_1$, further indicated by the high R$^2$ scores between the reconstructed and known points.}
\end{figure}

The presented implementation carries also reconstruction of scores $\mathbf{t}$ for a given $\mathbf{y}$ by the method \texttt{inverse()}. When only few ECA components cover most of the target variance, this procedure reduces the dimensionality of the inverse problem notably, potentially allowing for an approximate solution for it. Having multiple target values (dimensionality of $\mathbf{y}$) that provide independent information can be beneficial for this procedure, as is the case of this example. We generated a vector-valued $\mathbf{y}(\mathbf{x})$ by adding components to the function of Eqs. (\ref{eq:rudimentary1}) and (\ref{eq:rudimentary2}) depending only on $t_1$ as follows
\begin{equation}
\mathbf{y}(\mathbf{x}) 
= \begin{bmatrix}
(\mathbf{v}\cdot\mathbf{x})^3\\
\sin(0.2(\mathbf{v}\cdot\mathbf{x}))\\
\cos(0.2(\mathbf{v}\cdot\mathbf{x}))\\
\tanh(0.2(\mathbf{v}\cdot\mathbf{x}))
\end{bmatrix}
.\label{eq:rudimentary3}
\end{equation}
First, we z-score standardized the four features of the generated $\mathbf{y}$ data, and applied the same procedures as in the one-dimensional case to obtain emulators for each $d$ with R$^2$ score for the test set always better than 0.98. Then, we ran one dimensional ECA fit and projected all the ECA fit data points onto the obtained $\mathbf{v}_1$ providing us with scores $t_1$. Finally, we performed inverse ECA fit to reconstruct the $t_1$ scores of each data point. Both implementations reconstructed the scores well, as shown in Figure \ref{fig:inverse} for $d=2$ and $d=512$. The R$^2$ scores for $t_1$ between the projection and reconstruction are at least 0.98 for all input dimensionalities.

We observed constant inverse fit times regardless of $d$ for both implementations. However, the original implementation was faster by one order of magnitude. This is in agreement with the results in Figure \ref{fig:forward}c as the inverse problem is low dimensional and, therefore, does not benefit from the optimization done by PyTorch. Regardless, we consider the presented implementation relatively fast at around 0.07~s {\it per} reconstructed data point on the aforementioned cluster computer.

\section*{Scientific Example}
Next, we turn to the analysis of simulated N K-edge near-edge X-ray absorption fine structure (NEXAFS) of aqueous triglycine \cite{Eronen2024}. This data consists of local structural descriptors (local version of the many-body tensor representation, LMBTR \cite{Huo2022}, implemented in \cite{Himanen2020}) and related core-level spectra, available at \cite{triglydata}.  The descriptors are built from atomistic structures, consisting of the single triglycine molecule and the water molecules within 3\,\AA{} from it. The spectral data are given as intensities of the three regions of interest (ROIs; I--III). The purpose of the analysis is to discover spectrally dominant structural modifications, {\it i.e.} finding structural degrees the change of which results in the most change in the spectra. ECA was applied as in the original paper: the components were obtained with an ECA fit data set of 3000 data points and the results were validated with ECA validation data set of 2999 data points. The original emulator was trained on 23996 data points, and for this work it was simply converted (by replicating its architecture and weights) to PyTorch. The scripts of our implementation for this data are provided with the codes.
\par
\begin{figure*}
    \includegraphics[width=\linewidth]{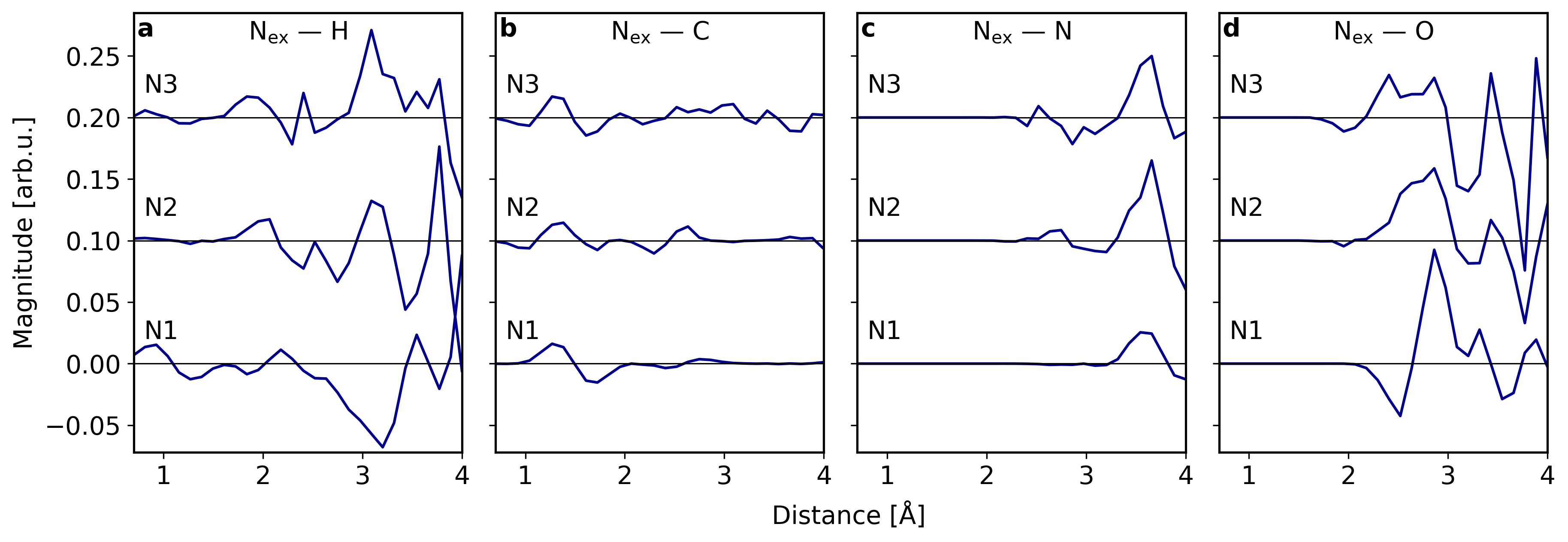}
    \caption{ The interatomic distance distribution part of the first component vector for the example of aqueous triglycine. The results shown here are in agreement with the original study in \cite{Eronen2024}. A non-zero value means that the studied spectral regions are sensitive to changes in this structural feature.}
    \label{fig:ECAcomps}
\end{figure*}

ECA covered variance up to rank five is given in Table~\ref{tab:ECAvariances} and closely resembles the original work. The LMBTR descriptor allows for an interpretation of the deduced basis vectors in terms of the structural features. Figure~\ref{fig:ECAcomps} presents the first basis vector which describes the interatomic distance distributions (from each active nitrogen site N$_\mathrm{ex}$). These results are also in agreement with the previous ones. Moreover, whereas the principal component analysis (PCA) of the structural descriptor shows two distinct regions, we showed that this information is lost in ROIs of NEXAFS, and therefore can not be recovered by X-ray absorption experiments. The results, available in the scripts provided with this work, were the same with the presented implementation.

\begin{table}[h]
\centering
\caption{The covered spectral ROI variances as a function of the rank of the ECA.}\label{tab:ECAvariances}
\begin{tabular}{ccc}
Rank & ECA & ECA \\
     & fit & validation \\
\hline
1   & 0.535 & 0.516 \\
2   & 0.747 & 0.710  \\ 
3   & 0.850 & 0.810 \\
4   & 0.865 & 0.811 \\
5   & 0.874 & 0.817 \\
\hline
\hline
\end{tabular}
\end{table}

\section{Discussion}
We find the results of the basis vector optimization to sometimes be dependent on the initial guess given for the vectors. Therefore optimization for each component vector should be carried out several ($\sim$10--100) times after which the best outcome is to be selected. This behaviour is understood by the loss function of the respective minimization problem being wiggly. Our selection of stochastic mini-batch optimization, although a viable means to reach global minimum, does not guarantee that this minimum is found from a random initial guess. The behaviour can likely be made more robust by choice of batch size, learning rate etc. We note that the solution of the inverse problem via inverse transformation did not seem to suffer from the aforementioned instability. 
\par
We discovered that the Adam optimizer works faster than the stochastic gradient descent (SGD) algorithm and reduces the dependence of the results from the initial guess (\textit{e.g.}, the random seed).  Utilization of the Adam optimizer requires less parameters to be specified, and in our tests with spectroscopic data we even found the default values to work well. We also discovered that for the test cases the learning rate for the inverse transformation must be rather large ($\ge0.05$). However, for robustness we would recommended to finely tune the optimizer hyperparameters as well as to check the reproducibility of the fit with, \textit{e.g.}, different initial guesses.

\par

When applied to X-ray spectroscopic data \cite{Niskanen2022,Vladyka2023,Eronen2024,Eronen2024b}, ECA deduced a subspace in structural data space which defines the majority of the spectral variations, {\it i.e.} represents the sensitivity of the spectral changes to the corresponding structural ones. The next step of the analysis is an interpretation of the obtained vectors as meaningful physical quantities. The example shown in this work is not a complete demonstration of the power the ECA -- other applications are possible and likely to arise in further studies. For example, these could include visualisation and analysis of the possible clustering in the data. Moreover, we have demonstrated that the magnitude of the first ECA component vector may be useful as a feature selection tool \cite{Eronen2024}.

Machine learning models, more specifically simple neural networks have been previously utilized inside variations of the PPR \cite{Lingjaerde1999,Barron1988,Intrator1993,Zhao1996}, but the applied nonlinear projection complicates the interpretation of the results which is essential for further analysis. ECA offers a simple interpretation of the relevant structural features through its component vectors: only a single score is required for each feature (per ECA component). This is further aided by our observation that often just rank one or two are needed to cover most of the target variance.

Although the schematic in Figure~\ref{fig:schematic} implies that the emulator is trained and the ECA is fitted with separate data sets, we have found no drawbacks when using the same data for both. After all, the emulator and the ECA can be considered as a single machine learning system.

\par


\section{Conclusions}
The implementation of the ECA algorithm as a Python class makes the method easily applicable for dimensionality reduction in inverse problems for which paired data of input and output values, and a computationally fast emulator, are available. The original application \cite{Niskanen2022} of the method targeted structural interpretation of X-ray spectra by a dimensionality reduction, which may be drastic \cite{Vladyka2023,Eronen2024,Eronen2024b}.  With the presented PyTorch implementation a robust and fairly fast ECA decomposition is available to any scientific Python user.

\section{Author Contributions}
A.V. programming, writing the manuscript. E.A.E. testing the software, writing the manuscript. J.N. research design, writing the manuscript, funding.

\section{Acknowledgements}
Academy of Finland is acknowledged for funding via project 331234. The authors acknowledge CSC – IT Center for Science, Finland, and the FGCI -- Finnish Grid and Cloud Infrastructure for computational resources. E.A.E. acknowledges Jenny and Antti Wihuri Foundation for funding. 

\section{Data availability}
The PyTorch implementation of the ECA and the relevant scripts are available in Zenodo: \href{https://doi.org/10.5281/zenodo.10405685}{10.5281/zenodo.10405685}. The data for the example of aqueous Triglycine N K-edge NEXAFS is available at \cite{triglydata}.

\bibliographystyle{unsrtnat}
\bibliography{references_revised}

\onecolumngrid
\newcommand{\indt}[1][1]{\parshape 1 #1\parindent \dimexpr\linewidth-#1\parindent\relax\noindent}
\newpage
\appendix
\section*{Reference manual}    

\noindent\textit{class} eca.ECA(emulator)

\noindent\textbf{\uline{Parameters}}:

\noindent\textbf{emulator}: \textit{torch.nn.Module}

Trained PyTorch-based emulator.
\newline\newline
\noindent\textbf{\uline{Default options}}:

\noindent DEFAULT\_LR = 1E-3

Learning rate. 

\noindent DEFAULT\_TOL = 1E-4

\indt Tolerance level. Optimization stops when the change of the calculated error (`loss') is smaller than tolerance,  \textit{i.e.} no more improvement in fit is observed.

\noindent DEFAULT\_EPOCHS = 10000

Maximum number of epochs for optimization. One epoch is an optimization step done for all available $\mathbf{X}$ data.

\noindent DEFAULT\_BETAS = (0.9, 0.999)

\indt Beta parameters for the Adam optimizer. These parameters are used for computing running averages of a gradient and its square (`momentum') \cite{Adam}. See help(torch.optim.Adam) for more details.

\noindent DEFAULT\_BATCH\_SIZE = 200

Size of mini-batches for optimization.

\noindent DEFAULT\_LR\_INV = 5E-2

Learning rate for inverse transformation.

\noindent DEFAULT\_TOL\_INV = 1E-4

Tolerance level for inverse transformation.

\noindent DEFAULT\_EPOCHS\_INV = 1000

Maximum number of epochs for inverse transformation. \\[0em]

\noindent\textbf{\uline{Attributes:}} (acquired during fit()):

\noindent\textbf{V}: \textit{torch.Tensor} -- calculated ECA components.

\noindent\textbf{y\_var}: \textit{list[float]} -- covered variance of $\mathbf{Y}$ data (cumulative).

\noindent\textbf{x\_var}: \textit{list[float]} -- covered variance of $\mathbf{X}$ data (cumulative).\\[0em]

\noindent\textbf{\uline{Methods}}:

\noindent\textbf{fit}(x, y, n\_comp=3, options=None, keep=0, verbose=False)

Performs ECA decomposition.

\textbf{Parameters}:

\quad\textbf{x}: \textit{torch.Tensor}, \textbf{y}: \textit{torch.Tensor} -- $\mathbf{X}$ and $\mathbf{Y}$ data, respectively. 

\quad\textbf{n\_comp}: \textit{int} -- number of components to evaluate.

\quad\textbf{options}: \textit{dict}

\quad\quad\textbf{`lr'}, \textit{float} -- learning rate.

\quad\quad\textbf{`betas'}, \textit{tuple(float, float)} -- `betas' parameters of Adam optimizer. 

\quad\quad\textbf{`tol'}, \textit{float} -- tolerance level.

\quad\quad\textbf{`epochs'}, \textit{int} -- maximum number of epochs.

\quad\quad\textbf{`batch\_size'}, \textit{int} -- size of mini-batches.

\quad\quad\textbf{`seed'}, \textit{int} -- seed for random number generator (affects initial guess and mini-batch splits).

\indt[2]\textbf{keep}: \textit{int} -- to keep specified number of components for another fit. \textit{E.g.}, keep=1 keeps only the first component; keep=$-2$ removes two last components.

\quad\textbf{verbose}: \textit{bool} -- if True, displays progress of the fit (requires \texttt{tqdm} package).

\textbf{Returns}:

\quad\textbf{V}: \textit{torch.Tensor} -- calculated ECA components.

\indt[2]\textbf{y\_var}: \textit{list[float]} -- covered variance of $\mathbf{Y}$ data (cumulative).

\indt[2]\textbf{x\_var}: \textit{list[float]} -- covered variance of $\mathbf{X}$ data (cumulative).\\[0em]

%
%
\noindent\textbf{inverse}(y, n\_comp=3, options=None, return\_error=False)

Perform inverse transformation from the given $\mathbf{y}$ to $\mathbf{t}'$.

\textbf{Parameters}:

\quad\textbf{y}: \textit{torch.Tensor} -- Y data for inverse transformation.

\indt[2]\textbf{n\_comp}: \textit{int} -- number of $\mathbf{t}'$ score components to evaluate, \textit{i.e.} the number of  $\mathbf{V}$ components to use for the inverse transformation. (None: use all available).

\quad\textbf{options}: \textit{dict} -- options for inverse transformation: learning rate and maximum number of epochs. See \textbf{fit()} for details. 

\textbf{Returns}:

\quad\textbf{t}: \textit{torch.Tensor} -- $\mathbf{t}'$ scores for given $\mathbf{y}$.

\quad\textbf{err}: \textit{torch.Tensor} -- mean-squared error for each evaluated data point. \\[0em]

%
%
\noindent\textbf{reconstruct}(y, n\_comp=3, options=None, return\_error=False)

Performs inverse transformation followed by its expansion to $\mathbf{X}$ space. \textbf{reconstruct(y)} $\equiv$ \textbf{expand(inverse(y))}.

\textbf{Parameters}:

\quad Same as for \textbf{inverse()}.

\textbf{Returns}:

\quad\textbf{x\_prime}: \textit{torch.Tensor} -- $\mathbf{x}'$ data for given $\mathbf{y}$.

\quad\textbf{err}: \textit{torch.Tensor} -- mean-squared error for each evaluated data point.\\[0em]

%
%

\noindent\textbf{transform}(x, n\_comp=None) 

Performs transformation of $\mathbf{x}$ data to $\mathbf{V}$ space.

\textbf{Parameters}:

\quad\textbf{x}: \textit{torch.Tensor} -- given $\mathbf{x}$ data.

\quad\textbf{n\_comp}: \textit{int} -- number of  $\mathbf{V}$ components to use for the transformation (None: use all available).

\textbf{Returns}:

\quad\textbf{t}: \textit{torch.Tensor} -- $\mathbf{t}$ scores for given $\mathbf{x}$.\\[0em]

%
%
\noindent\textbf{expand}(t, n\_comp=None) 

Expands $\mathbf{t}$ scores to $\mathbf{X}$ space.

\textbf{Parameters}:

\quad\textbf{t}: \textit{torch.Tensor} -- given $\mathbf{t}$ scores.

\quad\textbf{n\_comp}: \textit{int} -- number of  $\mathbf{V}$ components to use for the transformation (None: use all available).

\textbf{Returns}:

\quad\textbf{x\_prime}: \textit{torch.Tensor} -- $\mathbf{x}'$ data.\\[0em]
%
%

\noindent\textbf{project}(x, n\_comp=None)

Projects given $\mathbf{x}$ data to ECA space using specified number of $\mathbf{V}$ components. \textbf{project(y)} $\equiv$ \textbf{expand(transform(y))}.

\textbf{Parameters}:

\quad\textbf{x}: \textit{torch.Tensor} -- given $\mathbf{x}$ data.

\quad\textbf{n\_comp}: \textit{int} -- number of  $\mathbf{V}$ components to use for the transformation (None: use all available).

\textbf{Returns}:

\quad\textbf{x\_prime}: \textit{torch.Tensor} -- $\mathbf{x}'$ projections.\\[0em]
%
%

\noindent\textbf{set\_seed}(seed)

\indt Set seed for random number generator. Affects the initial guess as well as the mini-batches. Can be set explicitly before running \textbf{fit()} or \textbf{inverse()}.

\textbf{Parameters}:

\quad\textbf{seed}: \textit{int} -- seed.\\[0em]

%
%
\noindent\textbf{test}(x, y)

\indt For given ($\mathbf{x}, \mathbf{y}$) pair, calculates the covered variance $\rho(\mathbf{y}_\mathrm{pred}, \mathbf{y})$ between $\mathbf{y}$ and $\mathbf{y}_\mathrm{pred}$, predicted from the projection $\mathbf{x}'$ of $\mathbf{x}$ on ECA $\mathbf{V}$ space (see Eq.~\ref{R2}).

\textbf{Parameters}:

\quad\textbf{x, y}: \textit{torch.Tensor} -- given $\mathbf{x}$ and $\mathbf{y}$ data.

\textbf{Returns}:

\quad\textbf{rho}: \textit{float} -- covered variance.\\[0em]

\noindent\textbf{r2loss}(y\_pred, y\_known)

\indt Calculates missing variance between known $\mathbf{y}_\mathrm{known}$ and predicted $\mathbf{y}_\mathrm{pred}$ values defined as $1-\rho(\mathbf{y}_\mathrm{pred}, \mathbf{y}_\mathrm{known})$ (see Eq.~\ref{R2}).

\textbf{Parameters}:

\quad\textbf{y\_pred, y\_known}: \textit{torch.Tensor} -- given predicted and known values.

\textbf{Returns}:

\quad\textbf{r2}: \textit{torch.Tensor} -- covered variance.

\section*{Sample code to run ECA decomposition and inverse}

\begin{lstlisting}
## import necessary modules to handle data
import numpy as np
from eca import ECA
## Load the data, for example:
data_x = np.loadtxt(...)
data_y = np.loadtxt(...)
# convert data to torch.tensor if necessary: 
data_x_t, data_y_t = map(lambda x: torch.tensor(x, dtype=torch.float32), [data_x, data_y]) 
## Initialize an instance of the ECA class with an emulator
# For pytorch, a torch.nn.Module object can be called as an emulator as well as its forward() method
eca = ECA(model)
# define options
options = {
    "tol"       : 1E-4, 
    "epochs"    : 10000, 
    "lr"        : 1E-3, 
    "batch_size": 200,
    "seed"      : 123
}
# Carry out ECA decomposition
vectors, y_variance, x_variance = eca.fit(data_x_t, data_y_t, n_comp=3, options=options)

## Example: recalculate the last component with another seed
# set a new seed 
eca.set_seed(42)
# keep the first two (already calculated) components and redo the third one with a new seed
vectors, y_variance, x_variance  = eca.fit(data_x_t, data_y_t, n_comp=3, keep=-1)  

# Calculate scores for data_x
t_k = eca.transform(data_x_t)
# Directly project data_x on ECA subspace
data_x_k = eca.project(data_x_t)
# Search iteratively for scores matching the spectral data points
options_inv = {
    "tol"     : 1E-4, 
    "epochs"  : 1000, 
    "lr"      : 0.05, 
}
t_k_prime = eca.inverse(data_y_t, options=options_inv)

# Search for approximate data_x_appr matching the spectral data points. 
data_x_appr = eca.reconstruct(data_y_t)

\end{lstlisting}

\par

\end{document}